%% file: mart4.tex
\documentclass[12pt]{article}
\usepackage[francais,english]{babel}
\usepackage{amsmath}
\usepackage{amsfonts}
\usepackage{amssymb}
\usepackage{amsthm}
\usepackage{graphics}
\usepackage{amscd}
\usepackage[all]{xy}
\usepackage[T1]{fontenc}
\usepackage{epsfig}
\usepackage{color}
\input{psfig}
\newcommand{\Z}{\mathbb{Z}}
\newcommand{\R}{\mathbb{R}}
\newcommand{\N}{\mathbb{N}}

\newcommand{\C}{\mathbb{C}}

\newcommand{\E}{\mathbb{E}}

\renewcommand{\P}{\mathbb{P}}

\renewcommand{\H}{\mathbb{H}}

\newcommand{\A}{\mathcal{A}}
\newcommand{\SLE}{\operatorname{SLE}}

\DeclareMathOperator{\cotan}{cotan}

\title{SLE and triangles}
\author{Julien Dub\'edat\footnote {Universit\'e Paris-Sud}}

\newtheorem{Prop}{Proposition}
\newtheorem{Lem}{Lemma}
\newtheorem{Cor}{Corollary}
\newtheorem{Conj}{Conjecture}
\begin{document}

\maketitle

\begin{abstract}
By analogy with Carleson's observation on Cardy's formula describing 
crossing probabilities for the scaling limit of
critical percolation, we exhibit 
``privileged geometries'' for Stochastic Loewner Evolutions with
various parameters, for which 
certain hitting distributions are uniformly distributed.
We then examine consequences for limiting probabilities
 of events concerning 
various critical plane discrete models.
\end{abstract}

\section {Introduction}
It had been conjectured that 
many critical two-dimensional models from 
statistical physics
are conformally invariant in the scaling limit;
for instance, percolation, Ising/Potts models, FK percolation
 or dimers. The Stochastic Loewner
Evolution (SLE) introduced by Oded Schramm in \cite {S0}
 is a one-parameter family of
random paths in simply connected planar domains.
These processes are 
the only possible 
candidates for conformally invariant continuous limits of
the aforementioned discrete models.
See \cite {RS01} for a discussion of explicit conjectures.

Cardy \cite {Ca} used conformal field theory techniques 
to predict an explicit formula (involving a hypergeometric
function) that
should describe crossing probabilities of conformal rectangles
for critical percolation as a function of the aspect ratio of the
rectangle.
 Carleson pointed out that 
Cardy's formula could be expressed in a much simpler way by choosing 
another geometric setup, specifically by mapping the rectangle onto an
equilateral triangle $ABC$.
 The formula can then be simply described by saying that 
the probability of a crossing (in the triangle) between 
$AC$ and $BX$ for $X \in [BC]$ is $BX/BC$.  
Smirnov \cite {Sm1} proved rigorously Cardy's formula
for critical site percolation on the triangular lattice
and his proof
uses the global geometry of the equilateral
triangle (more than the local geometry of the triangular
lattice).

In the present paper, we
show that each $\SLE_\kappa$ is in some sense naturally associated 
to some geometrical normalization in that the formulas corresponding
to Cardy's formula can again be expressed in a simple way. 
Combining this with the conjectures on
continuous limits of various discrete models, this yields precise
simple conjectures on some asymptotics for
these models in particular geometric setups. Just as 
percolation may be associated with equilateral triangles, it 
turns that, for instance, the critical 2d Ising
model (and the FK percolation with parameter $q=2$) seems 
 to be associated
with right-angled isosceles triangles (because  $\SLE_{
  {16}/3}$ hitting probabilities in such triangles
are ``uniform'').
 Other isosceles triangles correspond to FK percolation with
different values of the $q$ parameter. In particular, $q=3$
corresponds to isosceles triangle with angle $2\pi/3$. 
Similarly, double dimer-models or $q=4$ Potts models (conjectured to
correspond to $\kappa=4$) seem to be best expressed in
strips (i.e., domains like $\R \times [0,1]$),
 and half-strips (i.e., $[0, \infty) \times [0,1]$)
are a favorable geometry for uniform
spanning trees.

\noindent{\bf Acknowledgments.} I wish to thank Wendelin Werner for
his help and advice, as well as Richard Kenyon for useful insight on
domino tilings and FK percolation models. I also wish to thank the
referee for numerous corrections and comments.

\section{Chordal SLE}
We first briefly 
recall the definition of chordal SLE in the upper half-plane
$\H$ going from $0$ to $\infty$ (see for instance 
\cite {LSW1,RS01} for more details).
 For any $z\in\H$, $t\geq 0$, define $g_t (z)$ by
$g_0(z)=z$
and 
$$\partial_tg_t(z)=\frac{2}{g_t(z)-W_t}$$
where $(W_t/\sqrt\kappa, t \ge 0)$ is a standard Brownian motion on $\R$,
starting from $0$. Let $\tau_z$ be the first time of explosion of this
ODE. Define the hull $K_t$ as
$$K_t=\overline{\{z\in\H \ : \  \tau_z<t\}}$$
The family $(K_t)_{t\geq 0}$ is an increasing family of compact sets in
$\overline{\H}$; furthermore, $g_t$ is a conformal equivalence of
$\H\backslash K_t$ onto $\H$. It has been proved (\cite{RS01}, 
see \cite{LSW2} for the case $\kappa=8$) that
there exists a continuous process $(\gamma_t)_{t\geq 0}$ with values
in $\overline{\H}$ such that $\H\backslash K_t$ is the unbounded
connected component of $\H \setminus\gamma_{[0,t]}$, a.s. This process is the
trace of the SLE and it can recovered from $g_t$ (and therefore from 
$W_t$) by
$$\gamma_t=\lim_{z\in\H\rightarrow W_t}g_t^{-1}(z)$$
For any simply connected domain $D$ with two boundary points (or prime
ends) $a$ and
$b$, chordal $\SLE_\kappa$ in $D$ from $a$ to $b$ is defined as
$K_t^{(D,a,b)}=h^{-1}(K_t^{(\H,0,\infty)})$, where
$K_t^{(\H,0,\infty)}$ is as above, and $h$ is a conformal equivalence
of $(D,a,b)$ onto $(\H,0,\infty)$. This definition is unambiguous up
to a linear time change thanks to the scaling property of SLE in the
upper half-plane (inherited from the scaling property of the driving process
$W_t$).

\section{A normalization of SLE}

The construction of SLE relies on 
the conformal equivalence $g_t$ of
$\H\backslash K_t$ onto $\H$. As $\H$ has non-trivial conformal
automorphisms, one can choose other conformal mappings. The original
$g_t$ is natural as all points of the real line seen from
infinity play the same role (hence the driving
process $(W_t)$ is a Brownian
motion). Other normalizations, such as the one used in $\cite{LSW1}$
may prove useful for different points of view.  

A by-product of Smirnov's results (\cite{Sm1}) is the following: let
$\kappa=6$, and $F$ be the conformal mapping of $(\H,0,1,\infty)$ onto
an equilateral triangle $(T,a,b,c)$. Let $h_t$ be the conformal
automorphism of $(T,a,b,c)$ such that $h_t(F(W_t))=a$, $h_t(F(g_t(1)))=b$,
$h_t(c)=c$. Then, for any $z\in\H$, $h_t(F(g_t(z)))$ is a
local martingale. 
Our goal in this section is to find similar functions $F$ 
for other values of $\kappa$.

Recall the definitions and notations of section 2. For $t<\tau_1$,
consider the conformal mapping of $\H\backslash K_t$ onto $\H$ defined
as:
$$\tilde g_t(z)=\frac{g_t(z)-W_t}{g_t(1)-W_t}$$
so that $\tilde g_t(\infty)=\infty$, $\tilde g_t(1)=1$ and $\tilde
g_t(\gamma_t)=0$, where $\gamma_t$ is the SLE trace.

Notice that if $F$ is an holomorphic map $D\rightarrow\C$ and
$(Y_t)_{t\geq 0}$ is a $D$-valued semimartingale, then (the bivariate
real version of) It\^o's formula
yields:
$$dF(Y_t)=\frac{dF}{dz}dY_t+\frac 12\frac{d^2F}{dz^2}d\langle Y_t\rangle$$
where the quadratic covariation $\langle.,.\rangle$ for real
semimartingales is extended in a
$\C$-bilinear fashion to complex semimartingales:
$$\langle Y_1,Y_2\rangle=(\langle\Re Y_1,\Re Y_2\rangle-\langle\Im Y_1,\Im Y_2\rangle)+i(\langle\Re Y_1,\Im
Y_2\rangle+\langle\Im Y_1,\Re Y_2\rangle)$$
so that $d\langle C_t\rangle=0$ for an isotropic complex Brownian motion
$(C_t)$. The setup here is slightly different from conformal
martingales as described in \cite{RY}.

In the present case, one gets:
$$d\tilde g_t(z)=\left[\frac 2{\tilde g_t(z)}-2\tilde
  g_t(z)+\kappa(\tilde
  g_t(z)-1)\right]\frac{dt}{(g_t(1)-W_t)^2}+(\tilde
g_t(z)-1)\frac{dW_t}{g_t(1)-W_t}$$
For notational convenience, define $w_t=\tilde g_t(z)$. After
performing the time change 
$$u(t)=\int_0^t\frac{ds}{(g_s(1)-W_s)^2}$$
one gets the autonomous SDE:
$$dw_u=(w_u-1)\left [\kappa-\frac 2{w_u}(1+w_u)\right ]du+(w_u-1)d\tilde W_u$$
where $(\tilde W_u/\sqrt\kappa)_{u\geq 0}$ is a standard Brownian
motion.

Let us take a closer look at the time change. Let
$Y_t=g_t(1)-W_t$; then, $dY_t=-dW_t + 2 dt / Y_t$,
so that $(Y_t / \sqrt {\kappa})_{t\geq 0}$ 
is a Bessel process of dimension $(1+
4 / \kappa)$. For $\kappa\leq 4$, 
this dimension is not smaller than $2$, so
that $Y$ almost surely never vanishes (see e.g. 
\cite {RY}); moreover,
a.s., 
$$\int_0^\infty\frac{dt}{Y_t^2}=\infty$$ 
Indeed, 
let $T_n=\inf \{ t > 0 \ : \ Y_t = 2^n \}$. Then, the positive random variables
$(\int_{T_n}^{T_{n+1}}{dt}/{Y_t^2} , n \ge 1)$ are i.i.d. (using the Markov
and scaling properties of Bessel processes).
Hence:
$$\int_0^\infty\frac{dt}{Y_t^2}\geq\sum_{n=1}^\infty\int_{T_n}^{T_{n+1}}\frac{dt}{Y_t^2}=\infty\text{\
  \ a.s.}$$
So the time change is a.s. a bijection from
$\R_+$ onto $\R_+$ if
$\kappa\leq 4$. 

When $\kappa>4$, the dimension of the Bessel
process $Y$ is smaller than $2$, so that $\tau_1<\infty$ almost 
surely. In
this case, using a similar argument with the stopping times $T_n$
for $n < 0$,
one sees that
$$\int_0^{\tau_1}\frac{dt}{Y_t^2}=\infty$$ 
Hence, if $\kappa>4$, the time change is a.s. a bijection
$[0,\tau_1)\rightarrow\R_+$.

We conclude that for all $\kappa >0$,
the stochastic flow $(\tilde g_u)_{u\geq 0}$ does
almost surely not explode in finite time.

We now look for holomorphic functions $F$ such that $(F(w_u))_{u\geq
  0}$ are local martingales. As before, one gets:
$$dF(w_u)=\left [F'(w_u)(\kappa-\frac 2{w_u}(1+w_u))+\frac \kappa 2
  F''(w_u)(w_u-1)\right](w_u-1)du+F'(w_u)(w_u-1)d\tilde W_u$$
Hence we have to find 
 holomorphic functions defined on $\H$
satisfying the following equation:
$$F'(w)\left[\kappa-\frac 2w(1+w)\right]+F''(w)\frac\kappa 2(w-1)=0$$
The solutions are such that 
$$F'(w)\propto
w^{\alpha-1}(w-1)^{\beta-1},$$ where
$$\left\{\begin{array}{lll}\alpha&=&1-\frac 4\kappa\\ \beta&=&\frac
    8\kappa-1\end{array}\right.$$
For $\kappa=4$, $F(w)=\log(w)$ is a solution.

\section{Privileged geometries}

In this section we attempt to identify the holomorphic
 map $F$ depending on the value of
the $\kappa$ parameter.

\begin{itemize}
\item {\bf Case $4<\kappa<8$}

Using the Schwarz-Christoffel formula \cite{Ahl}, one can identify $F$ as the
conformal equivalence of $(\H,0,1,\infty)$ onto an isosceles triangle
$(T_\kappa,a,b,c)$ with angles $\hat a=\hat c=\alpha\pi=(1-\frac 4\kappa)\pi$
 and $\hat b=\beta\pi=(\frac 8\kappa-1)\pi$. Special triangles turn
 out to correspond to special values of $\kappa$. Thus, for
 $\kappa=6$, one gets an equilateral triangle, as was foreseeable from
 Smirnov's work (\cite{Sm1}). For $\kappa=\frac{16}3$, a value conjectured
 to correspond to FK percolation with $q=2$ and to
 the Ising model, one
 gets an isorectangle triangle.

Since $F(\H)$ is bounded,  the local martingales
$F(\tilde g_{t\wedge\tau_1}(z))$ are bounded  
(complex-valued) martingales, so that one can
apply the optional stopping theorem. We therefore  study what
happens at the stopping time $\tau_{1,z}=\tau_1\wedge\tau_z$. There
are three possible cases, each having 
positive probability: $\tau_1<\tau_z$,
$\tau_1=\tau_z$ and $\tau_1>\tau_z$. Clearly, 
$\lim_{t\nearrow\tau_z}(g_t(z)-W_t)=0$, and on the other
hand $(g_t(z)-W_t)$ is bounded
away from zero if $t$ stays bounded away from $\tau_z$. Recall that
$$\tilde g_t(z)=\frac{g_t(z)-W_t}{g_t(1)-W_t}$$
So, as $t\nearrow\tau_{1,z}$,  $\tilde g_t(z)\rightarrow \infty$ if
$\tau_1<\tau_z$ and $\tilde g_t(z)\rightarrow 0$ if
$\tau_z<\tau_1$. 
In the case $\tau_1=\tau_z=\tau$, the points $1$ and $z$ are disconnected
at the same moment, with $\gamma_{\tau}\in\partial\H$. As
$t\nearrow\tau$, the harmonic measure of $(-\infty,0)$ seen from
$z$ tends to $0$; indeed, to reach $(-\infty,0)$, a 
Brownian motion starting from $z$ has to
go through the straits $[\gamma_t,\gamma_\tau]$ the width of which
tends to zero. At the same time, the harmonic measures of $(0,1)$ and
$(1,\infty)$ seen from $z$ stay bounded away from $0$. This implies that
$\tilde g_t(z)$ tends to $1$, as is easily seen by mapping $\H$ to strips.

Now one can apply the optional stopping theorem to the martingales
$F(\tilde g_{t\wedge\tau_{1,z}}(z))$. The mapping $F$ has a
continuous extension to $\overline\H$, hence:
$$F(z)=F(0)\P(\tau_z<\tau_1)+F(1)\P(\tau_z=\tau_1)+F(\infty)\P(\tau_z>\tau_1)$$
Thus:
\begin {Prop}\ \\
 The barycentric coordinates of $w=F(z)$ in the triangle
$T_\kappa$ are the probabilities of the events $\tau_z<\tau_1$,
$\tau_z=\tau_1$, $\tau_z>\tau_1$. 
\end {Prop}

Define $T^0=\{w\in T_\kappa \ : \ \tau_z<\tau_1\}$, $T^1=\{w\in
T_\kappa \ : \ \tau_z=\tau_1\}$, 
$T^\infty=\{w\in T_\kappa \ : \ \tau_z>\tau_1\}$,
which is a random partition of $T_\kappa$. These three sets are
a.s.\ borelian; indeed, $T^\infty=F(\H\backslash K_{\tau_1})$ is a.s. open, and
$T^0=\bigcup_{t<\tau_1}K_t$ is a.s. an $F_\sigma$ borelian. 
The integral of the above formula with respect to the Lebesgue measure
on $T_\kappa$ yields:

\begin{Cor}\ \\ The following relation holds:
$$\E(\A(T^0))= \E(\A(T^1))=\E(\A(T^\infty))=\frac{\A(T_\kappa)}3$$ 
where $\A$ designates the area.
\end{Cor}

\begin{figure}[htbp]
\begin{center}
\input{triang1.pstex_t}
\end{center}
\caption{The random partition}
\end{figure}
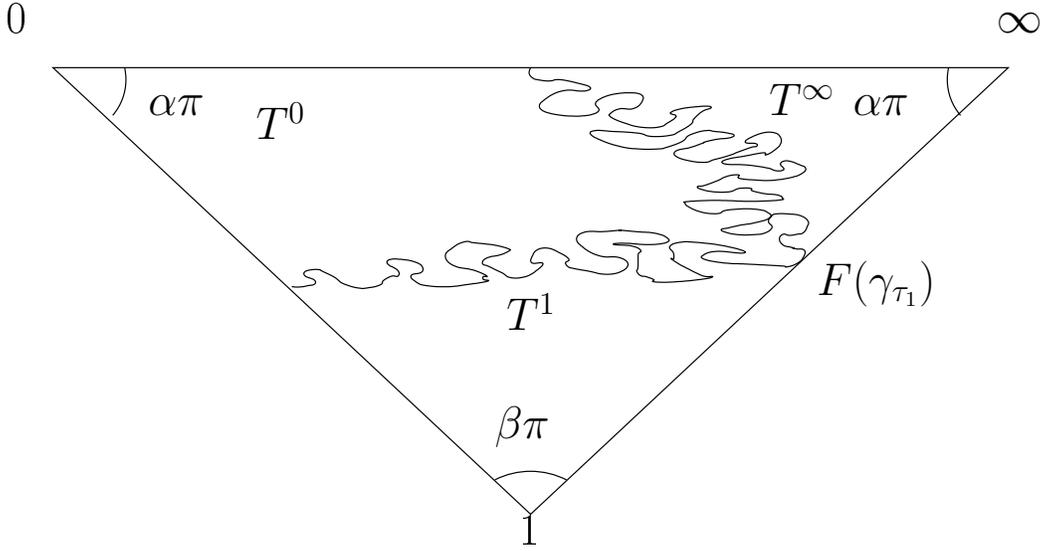

Another easy consequence is a Cardy's formula for SLE.
\begin{Cor}[Cardy's Formula]\ \\
Let $\gamma$ be the trace of a chordal $\SLE_\kappa$ going from $a$ to
$c$ in the isosceles triangle $T_\kappa$, $4<\kappa<8$. Let $\tau$ be the first time
$\gamma$ hits $(b,c)$. Then $\gamma_\tau$ has uniform distribution on $(b,c)$.
\end{Cor}
One can translate this result on the usual half-plane setup.
\begin{Cor}
Let $\gamma$ be the trace of a chordal $\SLE_\kappa$ going from $0$ to
$\infty$ in the half-plane, and $\gamma_{\tau_1}$ be the first hit of
the half-line $[1,\infty)$ by $\gamma$. Then, if $4<\kappa<8$, 
the 
law of $1/ 
\gamma_{\tau_1}$ is that of the beta distribution  $B(1-\frac 4\kappa,\frac 8\kappa-1)$. 
\end{Cor}
It is easy to see that, the law of
$\gamma_{\tau_1}$ converges weakly to $\delta_1$ when $\kappa \nearrow 8$.
 This is not
surprising as for  $\kappa\geq 8$, the
SLE trace $\gamma$ is a.s. a Peano curve, and $\gamma_{\tau_1}=1$ a.s.

\item {\bf Case $\kappa=4$}

In this case, $F(w)=\log(w)$ is a solution. One can choose a
determination of the logarithm such that $\Im(\log(\H))=(0,\pi)$. Then
$\Im(\log(\tilde g_t(z)))=\arg(\tilde g_t(z))$ is a bounded local
martingale. Let $\H_r$ (resp. $\H_l$) be the
points in $\H$ left on the right (resp. on the left) by the SLE trace
(a precise definition is to be found in \cite{S1}). If $z\in\H_l$, the
harmonic measure of $g_t^{-1}((W_t,\infty))$ seen from $z$ in
$\H\backslash\gamma_{[0,t]}$ tends to $0$ as
$t\rightarrow\infty=\tau_z$. This implies that the argument of $\tilde
g_t(z)$ tends to $\pi$. For $z\in\H_r$, an argument similar to the case
$4<\kappa<8$ shows that $\tilde g_t(z)\rightarrow 1$. 
Hence, applying the optional stopping theorem to the bounded
martingale $\arg(\tilde g_t(z))$, one gets:
$$\arg(z)=0 \times\P(z\in\H_r)+\pi\P(z\in\H_l)$$
or $\P(z\in\H_l)={\arg(z)} / \pi$, in accordance with \cite{S1}.

\begin{figure}[htbp]
\begin{center}
\input{triang4.pstex_t}
\end{center}
\caption{$F(\H)$, case $\kappa=4$: slit}
\end{figure}
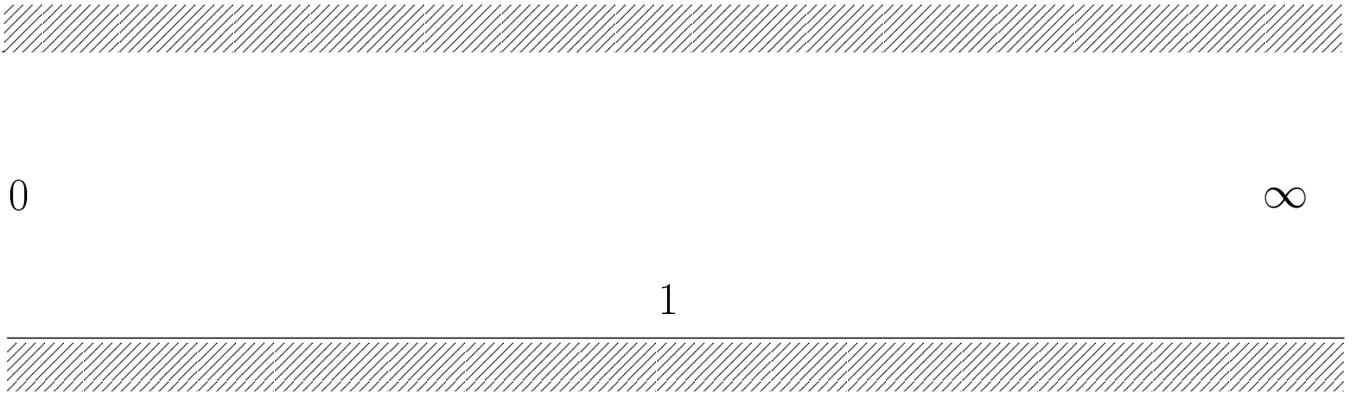

\item {\bf Case $\kappa=8$}

Let $F(z)=\int w^{-\frac 12}(w-1)^{-1}dw$; $F$ maps $(\H,0,1,\infty)$
onto a half-strip $(D,a,\infty,b)$. One may choose $F$ so that
$F(\H)=\{z \ : \ 
0<\Re z<1, \Im z>0\}$. Then $F(\infty)=0$ and $F(0)=1$. Moreover, $\Re F(\tilde g_t(z))$ is a bounded martingale. In
the case $\kappa\geq 8$, it is known that $\tau_1<\infty$,
$\tau_z<\infty$, and $\tau_1\neq\tau_z$ a.s. if $z\neq 1$ (see \cite{RS01}). Hence,
if $\tau=\tau_1\wedge \tau_z$, $\tilde g_\tau(z)$ equals
 $0$ or $\infty$, depending on whether $\tau_z<\tau_1$ or
 $\tau_z>\tau_1$. Applying the optional stopping theorem to the
 bounded martingale $\Re F(\tilde g_t(z))$, one gets:
$$\P(\tau_z<\tau_1)=\Re F(z)$$

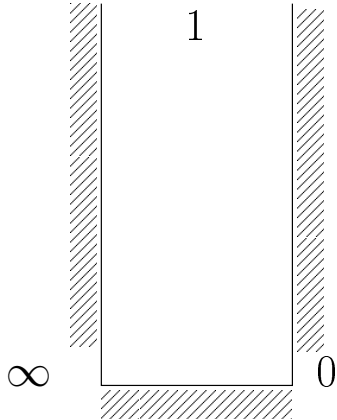
\begin{figure}[htbp]
\begin{center}
\input{triang5.pstex_t}
\end{center}
\caption{$F(\H)$, case $\kappa=8$: half-strip}
\end{figure}

\item{\bf Case $\kappa>8$}

In this case, one can choose $F$ so that it maps $(\H,0,1,\infty)$ onto
$(D,1,\infty,0)$ where 
$$D=\left\{z \ : \ \Im z>0, 0<\arg(z)<\left (1-\frac 4\kappa\right )\pi, \frac
4\kappa\pi<\arg(z-1)<\pi\right\}$$
Then $F(\H)$ is not bounded in any direction, preventing us from using the
optional stopping theorem.

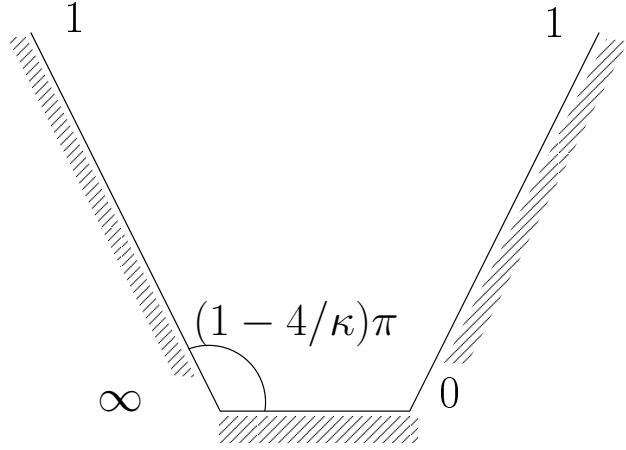
\begin{figure}[htbp]
\begin{center}
\input{triang3.pstex_t}
\end{center}
\caption{$F(\H)$, case $\kappa>8$}
\end{figure}

\item{\bf Case $\kappa<4$}

If $\kappa\geq 8/ 3$, one can choose $F$ so that it maps $(\H,0,1,\infty)$ onto
$(D,\infty,0,\infty)$, where
$$D=\left\{z \ : \ \Im z<1, -\left (\frac 4\kappa-1\right )\pi<\arg(z)<\frac
  4\kappa\pi\right\}$$
For $\kappa= 8/ 3$, one gets a slit half-plane. For $\kappa<
8/3$, the map $F$ ceases to be univalent. 

\begin{figure}[htbp]
\begin{center}
\input{triang6.pstex_t}
\end{center}
\caption{$F(\H)$, case $\frac 83\leq\kappa\leq 4$}
\end{figure}
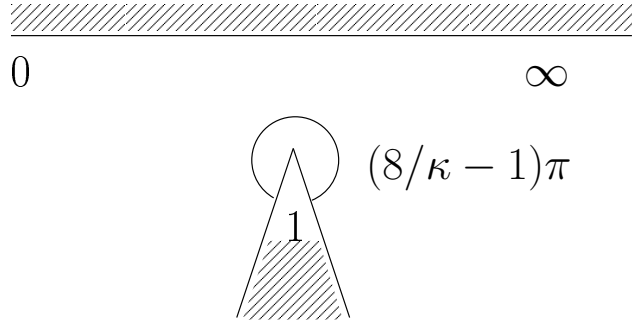

\end{itemize}

\section{Radial SLE}

Let $D$ be the unit disk. Radial SLE in $D$ starting from 1 is defined
by $g_0(z)=z$, $z\in D$ and the ODEs:
$$\partial_tg_t(z)=-g_t(z)\frac{g_t(z)+\xi(t)}{g_t(z)-\xi(t)}$$
where $\xi(t)=\exp(iW_t)$ and $W_t/\sqrt\kappa$ is a real standard
Brownian motion. The hulls $(K_t)$ and the trace $(\gamma_t)$ are
defined as in the chordal case (\cite{RS01}). 
Define $\tilde g_t(z)=g_t(z)\xi_t^{-1}$, so that $\tilde
g_t(0)=0$, $\tilde g(\gamma_t)=1$, where $(\gamma_t)$ is the SLE
trace. One may compute:
$$d\tilde g_t(z)=-\tilde g_t(z)\frac{\tilde g_t(z)+1}{\tilde g_t(z)-1}
dt+\tilde g_t(z)(-idW_t-\frac 12\kappa dt)$$
The above SDE is autonomous. 
As before, one looks for holomorphic functions $F$ such that $(F(\tilde
g_t(z)))_{t\geq 0}$ are local martingales. A sufficient condition is:
$$F'(z)\left(-z\frac{z+1}{z-1}-
\frac\kappa 2z\right)-\frac \kappa 2F''(z)z^2=0$$ 
i.e., 
$$\frac{F''(z)}{F'(z)}=\left(\frac 2\kappa-1\right)\frac 1z-\frac
4\kappa\frac 1{z-1}.$$
 Meromorphic solutions of this equation defined
on $D$ exist for $\kappa=2/n$, $n\in\N^*$. For $\kappa=2$,
$F(z)=(z-1)^{-1}$ is an (unbounded) solution.

\section{Related conjectures}

In this section we formulate various conjectures pertaining to
continuous limits of discrete critical models using the privileged
geometries for SLE described above.

\subsection{FK percolation in isosceles triangles}

For a survey of FK percolation, also called random-cluster model, see
\cite{G1}. We build on a conjecture stated in \cite{RS01}
(Conjecture 9.7),
according to which the discrete exploration process for critical FK percolation
with parameter $q$ converges weakly to the trace of $\SLE_\kappa$ for
$q\in (0,4)$, where the following relation holds:
$$\kappa=\frac{4\pi}{\cos^{-1}(-\sqrt q/2)}$$
Then the associated isosceles triangle $T_\kappa$ has angles $\hat
a=\hat c=\cos^{-1}(\sqrt q/2)$, $\hat b=\pi-2\hat a$. Let $\Gamma_n$
be a discrete approximation of the triangle $T_\kappa$ on the square
lattice with mesh $\frac 1n$; all vertices on the edges $(a,b]$ and
$[b,c)$ are identified. Let $\Gamma_n^\dagger$ be the dual graph. The
discrete exploration process $\beta$ runs between the opened connected component of
$(a,b]\cup [b,c)$ in $\Gamma_n$ and the closed connected component of
$(a,c)$ in $\Gamma_n^\dagger$. 

\begin{Conj}{Cardy's Formula}\ \\
Let $\tau$ be the first time $\beta$ hits $(b,c)$. Then, as $n$ tends
to infinity (i.e. as the mesh tends to zero), the law of $\beta_\tau$
converges weakly towards the uniform law on $(b,c)$.
\end{Conj}

\begin{figure}[htbp]
\begin{center}
\input{triang10.pstex_t}
\end{center}
\caption{Rectangle lattice, dual graph and associated isosceles triangle}\label{f10}
\end{figure}
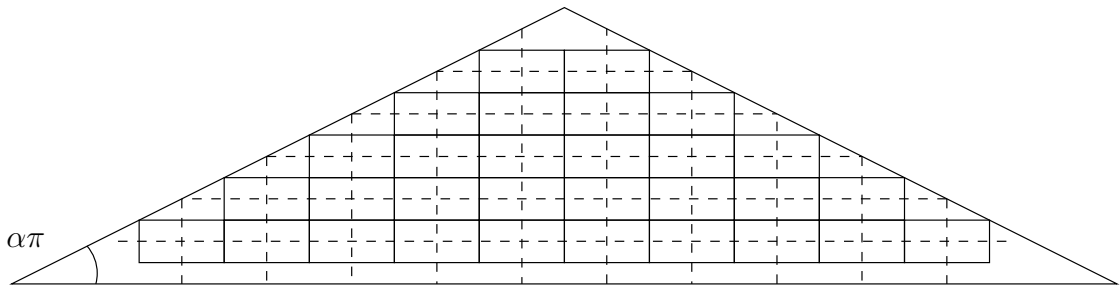

Kenyon
\cite {K3}
 has proposed an FK percolation model for any isoradial lattice,
in particular for any rectangular lattice. Let $\kappa$, $q$ and
$\alpha$ be as above, i.e. $4<\kappa<8$, $\frac{4\pi}\kappa=\cos^{-1}(-\sqrt q/2)$
and $\alpha=1-\frac 4\kappa$. Consider the rectangular lattice
$\Z\cos\alpha\pi+i\Z\sin\alpha\pi$. Then isosceles triangles homothetic
to $T_\kappa$ naturally fit in the lattice (see
figure \ref{f10}). Let $\Gamma=(V,E)$ be the finite graph  resulting from
the restriction of the lattice to a (large)
$T_\kappa$ triangle, with appropriate boundary conditions. A
configuration $\omega\in \{0,1\}^E$ of open edges has probability:
$$p_\Gamma(\omega)\propto q^{k(\omega)}\nu_h^{e_h(\omega)}\nu_v^{e_v(\omega)}$$
where $k(\omega)$ is the number of connected components in the
configuration, and $e_h$ (resp. $e_v$) is the number of open horizontal
(resp. vertical) edges. The weights $\nu_h$, $\nu_v$ are given by the
formulas:
\begin{align*}
\nu_v&=\sqrt q\frac{\sin (2\alpha^2\pi)}{\sin (\alpha(1-2\alpha)\pi)}\\
\nu_h&=\frac q{\nu_v}
\end{align*}

For this model, one may conjecture Cardy's formula as stated
above. Note that for $q=2$, $\kappa=\frac {16}3$, one retrieves the
usual critical FK percolation on the square lattice.

Let us now focus on the integral values of the $q$ parameter. It is known
that for these values there exists a stochastic coupling between FK
percolation and the Potts model (with parameter $q$) (see \cite{G1}).

\begin{itemize}
\item {$q=1$}\\
In this case FK percolation is simply percolation, $\kappa=6$, and the
privileged geometry is the equilateral triangle. This corresponds to
Carleson's observation on Cardy's formula.

\item{$q=2$}

\begin{figure}[htbp]
\begin{center}
\input{triang2.pstex_t}
\end{center}
\caption{Discrete exploration process for FK percolation ($q=2$, $\kappa=\frac{16}3$)}
\end{figure}
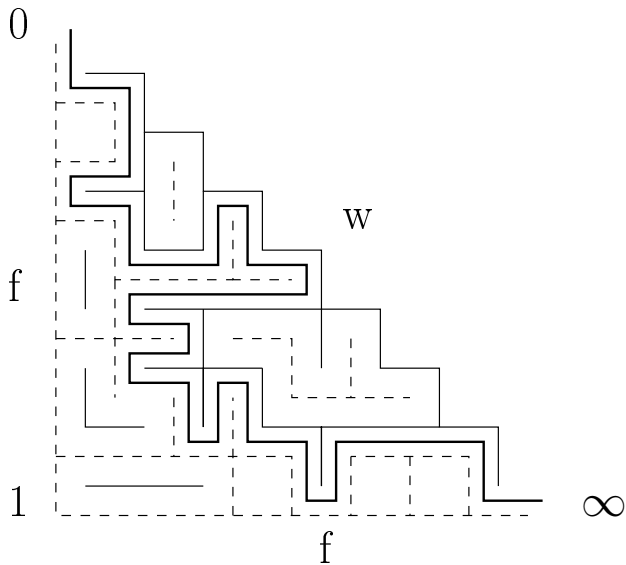
\noindent Here $\kappa=\frac{16}3$, and $T_\kappa$ is an isorectangle
triangle. As there is a stochastic coupling between FK percolation
with parameter $q=2$ and the Ising model (Potts model with $q=2$),
this suggests that the isorectangle triangle may be of some
significance for the Ising model. 

\item{$q=3$}\\
The corresponding geometry is the isosceles triangle $T_\frac{24}5$,
which has angles $\hat a=\hat c=\frac\pi 6$, $\hat b=\frac{2\pi}3$. The
possible relationship with the $q=3$ Potts model is not clear, as this
model is not
naturally associated with an exploration process.
\end{itemize}

\subsection{UST in half-strips}

It is proved in \cite{LSW2} that the scaling limit of the uniform spanning
trees (UST) Peano curve is the $\SLE_8$ chordal path. Let $R_{n,L}$ be
the square lattice $[0,n]\times[0,nL]$, with the following boundaries
conditions: the two horizontal arcs as well as the top one are wired,
and the bottom one is free. In fact, as we will consider the limit as $L$ goes to
infinity, one may as well consider that the top arc is free, which
makes the following lemma neater. We consider the uniform spanning tree in $R_{n,L}$. Let
$w$ be a point of the half-strip $\{z\ :\ 0<\Re z<1,\Im z>0\}$, and $w_n$ an
integral approximation of $nw$. Let
$a\in[0,n]$ be the unique triple point of the minimal subtree $T$
containing $(0,0)$, $(n,0)$ and $w_n$, and let $b$ be the other triple point
of the minimal subtree containing $(0,0)$, $(n,0)$, $w_n$ and
$(0,nL)$. One can formulate the following easy consequence of the
identification of the scaling limit of the UST:

\begin{Lem}\ \\
The following limits hold:
\begin{align*}
\lim_{n\rightarrow\infty}\lim_{L\rightarrow\infty}\P_{R_{n,L}}(b
\text{\rm\ belongs\ to\ the\ oriented\ arc\ }[0,a]\cup [a,w_n]\text{\ in\ }T)=\Re w\\
\lim_{L\rightarrow\infty}\lim_{n\rightarrow\infty}\P_{R_{n,L}}(b
\text{\rm\ belongs\ to\ the\ oriented\ arc\ }[0,a]\cup [a,w_n]\text{\ in\
}T)=\Re w
\end{align*}
\end{Lem}

Let us clarify the alternative (up to events of negligible
probability): either $b$ belongs to the (oriented) arc $[0,a]\cup
[a,w_n]$, or to the (oriented) arc $[w_n;a]\cup[a,1]$.
Recall that we have computed $\P(\tau_{F^{-1}(w)}>\tau_1)=\Re w$ for a chordal
$\SLE_8$ going from $0$ to $1$ in the half-strip (in accordance with
earlier conventions, subscripts refer to
points in the half-plane, not in the half-strip). As this path is
identified as the scaling limit of the UST Peano curve (start from
$0$ and go to $1$ with the UST rooted on the bottom always on your
right-hand), the event $\{\tau_1<\tau_{F^{-1}(w)}\}$ appears as a scaling limit
of an event involving only the subtree $T$. If one removes the arc
joining $a$ to $iLn$, $w_n$ is either on the left connected
component or on the right one depending on whether $w_n$ is
``visited'' by the exploration process before or after the top arc, up
to events of negligible probaility.

\begin{figure}[htbp]
\begin{center}
\input{triang7.pstex_t}
\end{center}
\caption{The alternative}
\end{figure}

In fact, one can prove the lemma without using the continuous
limit for UST. Indeed, let $w_n^\dagger$ be a point on the dual grid standing at
distance $\frac {\sqrt 2}2$ from $w_n$. Then, as $n$ tends to infinity,
\begin{align*}
\P_{R_{n,L}}(b
\text{\ belongs\ to\ the\ arc\ }[0,a]\cup [a,w_n]\text{\ in\ 
  }T)&\\
 -\P_{R_{n,L}^\dagger}(w_n^\dagger\text{\ is\ connected\ to\ the\ 
  right-hand\ boundary\ in\ the\ dual\ tree})&\rightarrow 0
\end{align*}

According to Wilson's algorithm \cite{W96}, the minimal subtree in the dual tree
connecting $w_n^\dagger$ to the boundary has the law of a loop-erased
random walk (LERW) stopped at its first hit of the boundary. The probability of
hitting the right-hand boundary or the left-hand boundary for a LERW
equals the corresponding probability for a simple random walk. The
continuous limit for a simple random walk with these boundary
conditions is a Brownian motion reflected on the bottom of the
half-strip; as the harmonic measure of the right-hand boundary of the
whole slit $\{0<\Re z<1\}$ seen
from $w_n^\dagger$ is $\Re w+o(1)$, this proves the lemma.

\subsection{Double domino tilings in plane strips}
For an early discussion of the double domino tiling model, see \cite{RHA}. It
is conjectured that the scaling limit of the path arising in this
model is the $\SLE_4$ trace (see \cite{RS01}, Problem 9.8). Building on
Kenyon's work \cite{K1,K2}, we show that the continuous
limit of a particular discrete event is 
compatible with the $\SLE_4$ conjecture.

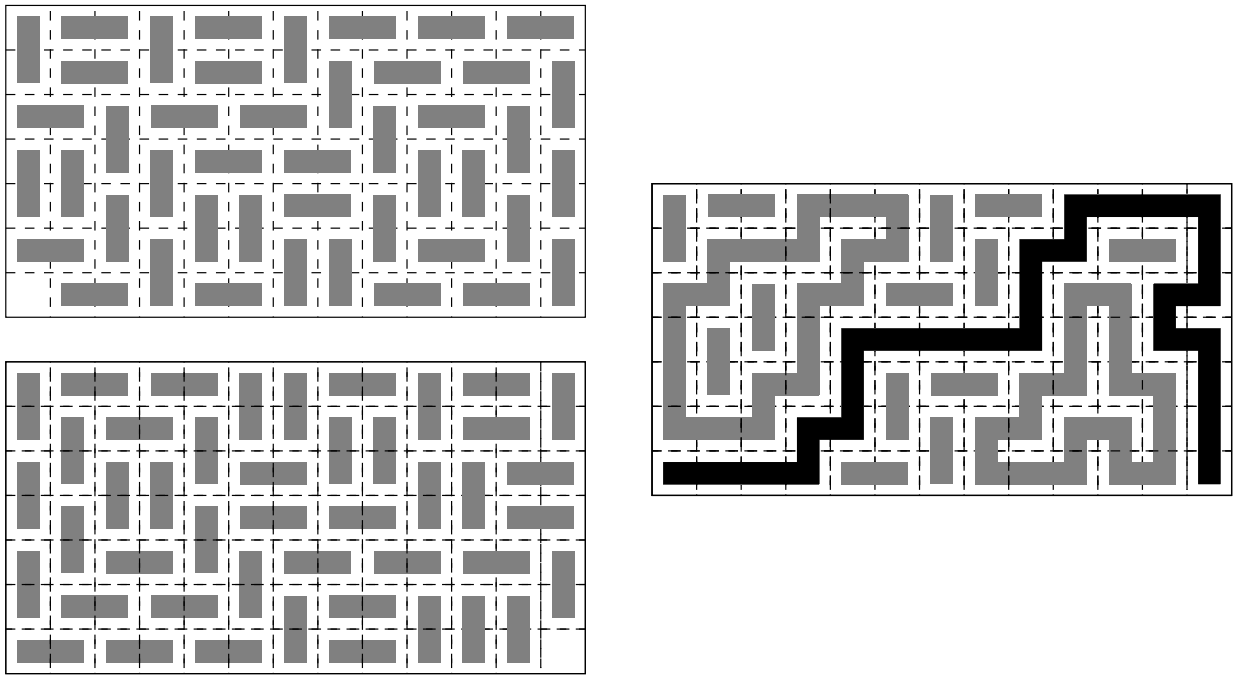
\begin{figure}[htbp]
\begin{center}
\input{triang9.pstex_t}
\end{center}
\caption{Double domino tilings and associated path}
\end{figure}

Consider the rectangle $R_{n,L}=[-nL,nL+1]\times [0,2n+1]$ (it is
important that the rectangle have odd length and width). Remove a unit
square at the corner $(-nL,0)$ or $(nL+1,0)$ to get two Temperleyan
polyominos (for general background on domino tilings, see \cite{K2}). Let $\gamma$ be the
random path going from $(-nL,0)$ to $(nL+1,0)$, arising from the superposed uniform domino tilings
on the two polyominos. Let $w$
be a point of the strip $\{z\ :\ 0<\Im z<1\}$, and $w_n$ an integral
approximation of $2nw$ in $R_{n,L}$.

\begin{Prop}\ \\
The following limit holds:
$$\lim_{L\rightarrow\infty}\lim_{n\rightarrow\infty}\P_{R_{n,L}}(w_n\text{\rm
  \ lies\ above\ }\gamma)=\Im z$$
\end{Prop}
\begin{proof}
We use a similar argument to the one given in \cite{K1}, 4.7. Let
$R_1$, $R_2$ be the two polyominos, and $h_1$, $h_2$ the height functions
associated with the two polyominos (these random integer-valued functions are defined up to
a constant). It is easily seen that one may choose $h_1$, $h_2$ so
that $h=h_1-h_2=0$ on the bottom side, and $h=4$ on the three other sides. Let
$x$ be an inner lattice point. Then:
$$\E(h(x))=4\P(x\text{ \ lies\ above\ }\gamma)$$
Indeed, condition on the union of the two dominos tilings. This union
consists of the path $\gamma$, doubled dominos and disjoint
cycles. Then $x$ is
separated from the bottom side by a certain number of closed cycles,
and possibly $\gamma$. Conditionally on the union, each closed cycle
accounts for $\pm 4$ with equal probability in
$h(x)$. Moreover, crossing $\gamma$ from below increases $h$ by
$4$. This yields the formula.

As $n$ goes to infinity, the average height functions converge to
harmonic functions (\cite{K2}, Theorem 23). Then take the limit as $L$
goes to infinity to conclude (one may map any finite rectangle $R_L$
to the whole slit, fixing a given point $x$; the boundary conditions
converge to the appropriate conditions, one concludes with Poisson's formula).
\end{proof}

\section{$\SLE(\kappa,\rho)$ processes and general triangles}

In this section we quickly discuss how any triangle may be associated
with a certain $\SLE$ process, in the same way as isoscele triangles
were associated with $\SLE_\kappa$ processes. 

\subsection{$\SLE(\kappa,\rho)$ processes}

Let us briefly describe $\SLE(\kappa,\rho)$ processes, defined in
\cite{LSW3}. Let $(W_t,O_t)_{t\geq 0}$ be a two-dimensional semimartingale satisfying
the following SDEs:
\begin{equation}\label{E11}
\left\{\begin{array}{l}dW_t=\sqrt\kappa dB_t+\frac\rho{W_t-O_t}dt\\
dO_t=\frac 2{O_t-W_t}dt\end{array}\right.
\end{equation}
where $B$ is a standard Brownian motion, as well as the inequality
$W_t\leq O_t$ valid for all positive times (the convention here
differs from the one in \cite{LSW3}). This process is well defined
for $\kappa>0$, $\rho>-2$. Indeed, one may consider
$Z_t=O_t-W_t$. The process $(Z_t/\sqrt\kappa)_{t\geq 0}$ is a Bessel process
in dimension $d=1+2\frac{\rho+2}\kappa$. Such processes are well
defined semimartingales if $d>1$, or $\rho>-2$ (see for instance \cite{RY}). Then
$O_t=2\int_0^t\frac{du}{Z_u}$ and $W_t=O_t-Z_t$. 

Hence one may define a $\SLE(\kappa,\rho)$ as a stochastic Loewner
chain the driving process of which has the law of the process $(W_t)$
defined above. The starting point (or rather state) of the process is
a couple $(w,o)$ with $w\leq o$, usually set to $(0,0^+)$. 
Then $O_t$ represents the image under the conformal mapping 
$g_t$ of the rightmost point of $\partial K_t\cup {O_0}$. Obviously,
for $\rho=0$, one recovers a standard $\SLE(\kappa)$ process.

\begin{Prop}
Let $(W_t,O_t)$ be the driving process of a $\SLE(\kappa,\rho)$
process starting from $(0,1)$, and $(g_t)$ be the associated conformal
equivalences. Let $z\in\H$. Then if $F$ is any analytic function on
$\H$, the complex-valued semimartingale
$$t\mapsto F\left(\frac{g_t(z)-W_t}{O_t-W_t}\right)$$
is a local martingale if and only if:
$$F'(z)\propto z^{-\frac 4\kappa}(1-z)^{2\frac{\rho-\kappa+4}\kappa}.$$
\end{Prop}
The proof is routine and is omitted. Once again, the conformal mapping $F$ may be
identified using the Schwarz-Christoffel formula.

\subsection{A particular case}

In \cite{S1}, Schramm derives expressions of the form
$$\P(z\in\H \text{\ lies\ to\ the\ left\ of\ }\gamma)=F_\kappa(\arg z)$$
where $\gamma$ is the trace of a $\SLE(\kappa)$ process, for
$\kappa\leq 4$. The function $F_\kappa$ involves hypergeometric functions,
and $F_\kappa(x)\propto x$ iff $\kappa=4$ (in this case $F\circ\arg$ is a
harmonic function). Now it is easily seen that for any $\kappa>0$,
$\rho>-2$, if $\delta$ designates the right boundary
of a $\SLE(\kappa,\rho)$ process starting from $(0,0^+)$, then a
simple consequence of scaling is the existence of a function
$F_{\kappa,\rho}$ such that:
$$\P(z\in\H \text{\ lies\ to\ the\ left\ of\ }\delta)=F_{\kappa,\rho}(\arg z)$$
Moreover, this function is not identically zero if $\rho\geq\kappa/2-2$.
This motivates the following result:
\begin{Prop}
Let $\kappa>4$, $\rho=\frac\kappa 2-2$. Then:
$$\P(z\in\H \text{\rm\ lies\ to\ the\ left\ of\ }\delta)=\arg z/\pi$$
\end{Prop}
\begin{proof}
Lying to the left of the right boundary of the hull is the same thing
as being absorbed if $\kappa>4$. Let $(W_t,O_t)$ be the driving
mechanism of the $\SLE(\kappa,\frac\kappa 2-2)$, and let
$z_t=g_t(z)$. Suppose for now that the starting state of the $\SLE$ is
$(W_0,O_0)=(0,1)$. Let $h:\H\rightarrow\C$ be a holomorphic
function. We have seen that a necessary and sufficient
condition for $h(\frac{z_t-W_t}{O_t-W_t})$ to be a ($\C$-valued) local
martingale is the holomorphic differential equation:
$$\frac {h''(z)}{h'(z)}=-\frac 4\kappa\frac
1z-2\frac{\rho-\kappa+4}\kappa\frac1{1-z}$$
or $h(z)\propto z^{-\frac
  4\kappa}(1-z)^{2\frac{\rho-\kappa+4}\kappa}$. In the case
$\rho=\kappa/2-2$, using the
Schwarz-Christoffel formula (see \cite{Ahl}), one sees that $h$ is (up
to a constant factor) the conformal equivalence between
$(\H,0,1,\infty)$ and $(D,0,1,\infty)$, where $D$ is the degenerate
triangle defined by:
$$D=\{z\in\H\ :\ \arg(z)\leq\pi(1-4/\kappa),\arg(z-1)\geq\pi(1-4/\kappa)\}$$
Let $\varphi(z)=\Re z-\cotan(\pi(1-4/\kappa))\Im z$. Then the image of
$D$ under this $\R$-linear form is $[0,1]$. Hence $\varphi\circ
h(\frac{z_t-W_t}{O_t-W_t})$ is a bounded martingale. Moreover,
standard convergence arguments imply that $\frac{z_t-W_t}{O_t-W_t}$
goes to $0$ in finite time if $z$ is absorbed and to 1 in infinite
time in the other case. A straightforward application of the optional
stopping theorem yields:
$$\P(z\in\H \text{\ lies\ to\ the\ right\
  of\ }\delta)=\varphi\left(\int_0^z w^{-\frac 4\kappa}(1-w)^{\frac 4\kappa-1}dw\right)/B(1-4/\kappa,4/\kappa)$$
Taking the asymptotics of this formula when $z=r\exp^{i\theta}$ goes to
infinity with constant argument (making use of $B(1-x,x)=\pi/\sin(\pi x)$), one finds that for a
$\SLE(\kappa,\frac\kappa 2-2)$ starting from $(0,0^+)$ :
$$\P(z\in\H \text{\ lies\ to\ the\ right\
  of\ }\delta)=1-\arg z/\pi$$
\end{proof}

In other words, $F_{\kappa,\kappa/2-2}=F_4$ for all $\kappa\geq
4$. This raises several questions, such as whether this still holds
for $\kappa<4$, or whether in full generality
$F_{\kappa,\rho}=F_{2\kappa/(\rho+2)}$, this last conjecture being
based on the dimension of the Bessel process $(O_t-W_t)$, where
$(W_t,O_t)$ designates the driving process of a $\SLE(\kappa,\rho)$ process.

-----------------------

Laboratoire de Math\'ematiques, B\^at. 425

Universit\'e Paris-Sud, F-91405 Orsay cedex, France

julien.dubedat@math.u-psud.fr

\end{document}

%% file: triang1.pstex_t
\begin{picture}(0,0)%
\special{psfile=triang1.pstex}%
\end{picture}%
\setlength{\unitlength}{3947sp}%
\begingroup\makeatletter\ifx\SetFigFont\undefined%
\gdef\SetFigFont#1#2#3#4#5{%
  \reset@font\fontsize{#1}{#2pt}%
  \fontfamily{#3}\fontseries{#4}\fontshape{#5}%
  \selectfont}%
\fi\endgroup%
\begin{picture}(6312,3429)(2701,-3361)
\put(2701,-136){\makebox(0,0)[lb]{\smash{\SetFigFont{17}{20.4}{\rmdefault}{\mddefault}{\updefault}{\color[rgb]{0,0,0}0}%
}}}
\put(8926,-136){\makebox(0,0)[lb]{\smash{\SetFigFont{17}{20.4}{\rmdefault}{\mddefault}{\updefault}{\color[rgb]{0,0,0}$\infty$}%
}}}
\put(4276,-811){\makebox(0,0)[lb]{\smash{\SetFigFont{17}{20.4}{\rmdefault}{\mddefault}{\updefault}{\color[rgb]{0,0,0}$T^0$}%
}}}
\put(7801,-1786){\makebox(0,0)[lb]{\smash{\SetFigFont{17}{20.4}{\rmdefault}{\mddefault}{\updefault}{\color[rgb]{0,0,0}$F(\gamma_{\tau_1})$}%
}}}
\put(5926,-3361){\makebox(0,0)[lb]{\smash{\SetFigFont{17}{20.4}{\rmdefault}{\mddefault}{\updefault}{\color[rgb]{0,0,0}1}%
}}}
\put(5776,-2686){\makebox(0,0)[lb]{\smash{\SetFigFont{17}{20.4}{\rmdefault}{\mddefault}{\updefault}{\color[rgb]{0,0,0}$\beta\pi$}%
}}}
\put(7501,-661){\makebox(0,0)[lb]{\smash{\SetFigFont{17}{20.4}{\rmdefault}{\mddefault}{\updefault}{\color[rgb]{0,0,0}$T^\infty$}%
}}}
\put(5851,-2011){\makebox(0,0)[lb]{\smash{\SetFigFont{17}{20.4}{\rmdefault}{\mddefault}{\updefault}{\color[rgb]{0,0,0}$T^1$}%
}}}
\put(3601,-661){\makebox(0,0)[lb]{\smash{\SetFigFont{17}{20.4}{\rmdefault}{\mddefault}{\updefault}{\color[rgb]{0,0,0}$\alpha\pi$}%
}}}
\put(8026,-661){\makebox(0,0)[lb]{\smash{\SetFigFont{17}{20.4}{\rmdefault}{\mddefault}{\updefault}{\color[rgb]{0,0,0}$\alpha\pi$}%
}}}
\end{picture}

%% file: triang4.pstex_t
\begin{picture}(0,0)%
\special{psfile=triang4.pstex}%
\end{picture}%
\setlength{\unitlength}{3947sp}%
\begingroup\makeatletter\ifx\SetFigFont\undefined%
\gdef\SetFigFont#1#2#3#4#5{%
  \reset@font\fontsize{#1}{#2pt}%
  \fontfamily{#3}\fontseries{#4}\fontshape{#5}%
  \selectfont}%
\fi\endgroup%
\begin{picture}(8454,2432)(1174,-1892)
\put(5311,-1411){\makebox(0,0)[lb]{\smash{\SetFigFont{17}{20.4}{\rmdefault}{\mddefault}{\updefault}{\color[rgb]{0,0,0}1}%
}}}
\put(1231,-751){\makebox(0,0)[lb]{\smash{\SetFigFont{17}{20.4}{\rmdefault}{\mddefault}{\updefault}{\color[rgb]{0,0,0}0}%
}}}
\put(9106,-736){\makebox(0,0)[lb]{\smash{\SetFigFont{17}{20.4}{\rmdefault}{\mddefault}{\updefault}{\color[rgb]{0,0,0}$\infty$}%
}}}
\end{picture}

%% file: triang5.pstex_t
\begin{picture}(0,0)%
\special{psfile=triang5.pstex}%
\end{picture}%
\setlength{\unitlength}{3947sp}%
\begingroup\makeatletter\ifx\SetFigFont\undefined%
\gdef\SetFigFont#1#2#3#4#5{%
  \reset@font\fontsize{#1}{#2pt}%
  \fontfamily{#3}\fontseries{#4}\fontshape{#5}%
  \selectfont}%
\fi\endgroup%
\begin{picture}(1996,2623)(3001,-2372)
\put(4951,-2161){\makebox(0,0)[lb]{\smash{\SetFigFont{17}{20.4}{\rmdefault}{\mddefault}{\updefault}{\color[rgb]{0,0,0}0}%
}}}
\put(4126, 14){\makebox(0,0)[lb]{\smash{\SetFigFont{17}{20.4}{\rmdefault}{\mddefault}{\updefault}{\color[rgb]{0,0,0}1}%
}}}
\put(3001,-2161){\makebox(0,0)[lb]{\smash{\SetFigFont{17}{20.4}{\rmdefault}{\mddefault}{\updefault}{\color[rgb]{0,0,0}$\infty$}%
}}}
\end{picture}

%% file: triang3.pstex_t
\begin{picture}(0,0)%
\special{psfile=triang3.pstex}%
\end{picture}%
\setlength{\unitlength}{3947sp}%
\begingroup\makeatletter\ifx\SetFigFont\undefined%
\gdef\SetFigFont#1#2#3#4#5{%
  \reset@font\fontsize{#1}{#2pt}%
  \fontfamily{#3}\fontseries{#4}\fontshape{#5}%
  \selectfont}%
\fi\endgroup%
\begin{picture}(3902,2761)(2250,-2342)
\put(2626,239){\makebox(0,0)[lb]{\smash{\SetFigFont{17}{20.4}{\rmdefault}{\mddefault}{\updefault}{\color[rgb]{0,0,0}1}%
}}}
\put(2836,-2146){\makebox(0,0)[lb]{\smash{\SetFigFont{17}{20.4}{\rmdefault}{\mddefault}{\updefault}{\color[rgb]{0,0,0}$\infty$}%
}}}
\put(3436,-1666){\makebox(0,0)[lb]{\smash{\SetFigFont{17}{20.4}{\rmdefault}{\mddefault}{\updefault}{\color[rgb]{0,0,0}$(1-4/\kappa)\pi$}%
}}}
\put(5641,209){\makebox(0,0)[lb]{\smash{\SetFigFont{17}{20.4}{\rmdefault}{\mddefault}{\updefault}{\color[rgb]{0,0,0}1}%
}}}
\put(4981,-2116){\makebox(0,0)[lb]{\smash{\SetFigFont{17}{20.4}{\rmdefault}{\mddefault}{\updefault}{\color[rgb]{0,0,0}0}%
}}}
\end{picture}

%% file: triang6.pstex_t
\begin{picture}(0,0)%
\special{psfile=triang6.pstex}%
\end{picture}%
\setlength{\unitlength}{3947sp}%
\begingroup\makeatletter\ifx\SetFigFont\undefined%
\gdef\SetFigFont#1#2#3#4#5{%
  \reset@font\fontsize{#1}{#2pt}%
  \fontfamily{#3}\fontseries{#4}\fontshape{#5}%
  \selectfont}%
\fi\endgroup%
\begin{picture}(3924,1982)(1789,-2147)
\put(1801,-661){\makebox(0,0)[lb]{\smash{\SetFigFont{17}{20.4}{\rmdefault}{\mddefault}{\updefault}{\color[rgb]{0,0,0}0}%
}}}
\put(5026,-661){\makebox(0,0)[lb]{\smash{\SetFigFont{17}{20.4}{\rmdefault}{\mddefault}{\updefault}{\color[rgb]{0,0,0}$\infty$}%
}}}
\put(3526,-1636){\makebox(0,0)[lb]{\smash{\SetFigFont{17}{20.4}{\rmdefault}{\mddefault}{\updefault}{\color[rgb]{0,0,0}1}%
}}}
\put(4036,-1246){\makebox(0,0)[lb]{\smash{\SetFigFont{17}{20.4}{\rmdefault}{\mddefault}{\updefault}{\color[rgb]{0,0,0}$(8/\kappa-1)\pi$}%
}}}
\end{picture}

%% file: triang10.pstex_t
\begin{picture}(0,0)%
\special{psfile=triang10.pstex}%
\end{picture}%
\setlength{\unitlength}{3947sp}%
\begingroup\makeatletter\ifx\SetFigFont\undefined%
\gdef\SetFigFont#1#2#3#4#5{%
  \reset@font\fontsize{#1}{#2pt}%
  \fontfamily{#3}\fontseries{#4}\fontshape{#5}%
  \selectfont}%
\fi\endgroup%
\begin{picture}(6976,1761)(-187,-1643)
\put(-187,-1395){\makebox(0,0)[lb]{\smash{\SetFigFont{12}{14.4}{\rmdefault}{\mddefault}{\updefault}{\color[rgb]{0,0,0}$\alpha\pi$}%
}}}
\end{picture}

%% file: triang2.pstex_t
\begin{picture}(0,0)%
\special{psfile=triang2.pstex}%
\end{picture}%
\setlength{\unitlength}{3947sp}%
\begingroup\makeatletter\ifx\SetFigFont\undefined%
\gdef\SetFigFont#1#2#3#4#5{%
  \reset@font\fontsize{#1}{#2pt}%
  \fontfamily{#3}\fontseries{#4}\fontshape{#5}%
  \selectfont}%
\fi\endgroup%
\begin{picture}(3600,3480)(4201,-2911)
\put(4201,389){\makebox(0,0)[lb]{\smash{\SetFigFont{17}{20.4}{\rmdefault}{\mddefault}{\updefault}{\color[rgb]{0,0,0}0}%
}}}
\put(4201,-2611){\makebox(0,0)[lb]{\smash{\SetFigFont{17}{20.4}{\rmdefault}{\mddefault}{\updefault}{\color[rgb]{0,0,0}1}%
}}}
\put(7801,-2611){\makebox(0,0)[lb]{\smash{\SetFigFont{17}{20.4}{\rmdefault}{\mddefault}{\updefault}{\color[rgb]{0,0,0}$\infty$}%
}}}
\put(6301,-811){\makebox(0,0)[lb]{\smash{\SetFigFont{17}{20.4}{\rmdefault}{\mddefault}{\updefault}{\color[rgb]{0,0,0}w}%
}}}
\put(4201,-1261){\makebox(0,0)[lb]{\smash{\SetFigFont{17}{20.4}{\rmdefault}{\mddefault}{\updefault}{\color[rgb]{0,0,0}f}%
}}}
\put(6151,-2911){\makebox(0,0)[lb]{\smash{\SetFigFont{17}{20.4}{\rmdefault}{\mddefault}{\updefault}{\color[rgb]{0,0,0}f}%
}}}
\end{picture}

%% file: triang7.pstex_t
\begin{picture}(0,0)%
\special{psfile=triang7.pstex}%
\end{picture}%
\setlength{\unitlength}{3947sp}%
\begingroup\makeatletter\ifx\SetFigFont\undefined%
\gdef\SetFigFont#1#2#3#4#5{%
  \reset@font\fontsize{#1}{#2pt}%
  \fontfamily{#3}\fontseries{#4}\fontshape{#5}%
  \selectfont}%
\fi\endgroup%
\begin{picture}(7035,4620)(3916,-4036)
\put(4906,-3856){\makebox(0,0)[lb]{\smash{\SetFigFont{17}{20.4}{\rmdefault}{\mddefault}{\updefault}{\color[rgb]{0,0,0}a}%
}}}
\put(3916,-3976){\makebox(0,0)[lb]{\smash{\SetFigFont{17}{20.4}{\rmdefault}{\mddefault}{\updefault}{\color[rgb]{0,0,0}0}%
}}}
\put(6766,-3931){\makebox(0,0)[lb]{\smash{\SetFigFont{17}{20.4}{\rmdefault}{\mddefault}{\updefault}{\color[rgb]{0,0,0}n}%
}}}
\put(5521,-2461){\makebox(0,0)[lb]{\smash{\SetFigFont{17}{20.4}{\rmdefault}{\mddefault}{\updefault}{\color[rgb]{0,0,0}$w_n$}%
}}}
\put(5461,-3301){\makebox(0,0)[lb]{\smash{\SetFigFont{17}{20.4}{\rmdefault}{\mddefault}{\updefault}{\color[rgb]{0,0,0}b}%
}}}
\put(4141,344){\makebox(0,0)[lb]{\smash{\SetFigFont{17}{20.4}{\rmdefault}{\mddefault}{\updefault}{\color[rgb]{0,0,0}(0,Ln)}%
}}}
\put(9091,-3841){\makebox(0,0)[lb]{\smash{\SetFigFont{17}{20.4}{\rmdefault}{\mddefault}{\updefault}{\color[rgb]{0,0,0}a}%
}}}
\put(8101,-3961){\makebox(0,0)[lb]{\smash{\SetFigFont{17}{20.4}{\rmdefault}{\mddefault}{\updefault}{\color[rgb]{0,0,0}0}%
}}}
\put(10951,-3916){\makebox(0,0)[lb]{\smash{\SetFigFont{17}{20.4}{\rmdefault}{\mddefault}{\updefault}{\color[rgb]{0,0,0}n}%
}}}
\put(8326,359){\makebox(0,0)[lb]{\smash{\SetFigFont{17}{20.4}{\rmdefault}{\mddefault}{\updefault}{\color[rgb]{0,0,0}(0,Ln)}%
}}}
\put(9181,-2746){\makebox(0,0)[lb]{\smash{\SetFigFont{17}{20.4}{\rmdefault}{\mddefault}{\updefault}{\color[rgb]{0,0,0}$w_n$}%
}}}
\put(9676,-2086){\makebox(0,0)[lb]{\smash{\SetFigFont{17}{20.4}{\rmdefault}{\mddefault}{\updefault}{\color[rgb]{0,0,0}b}%
}}}
\end{picture}

%% file: triang9.pstex_t
\begin{picture}(0,0)%
\special{psfile=triang9.pstex}%
\end{picture}%
\setlength{\unitlength}{3947sp}%
\begingroup\makeatletter\ifx\SetFigFont\undefined%
\gdef\SetFigFont#1#2#3#4#5{%
  \reset@font\fontsize{#1}{#2pt}%
  \fontfamily{#3}\fontseries{#4}\fontshape{#5}%
  \selectfont}%
\fi\endgroup%
\begin{picture}(7723,4224)(1789,-4573)
\end{picture}